\documentclass[letterpaper]{amsart}

\usepackage[T1]{fontenc}
\usepackage{color}
\usepackage{amsrefs}
\usepackage{hyperref}
\hypersetup{breaklinks=true}
\usepackage[utf8]{inputenc} 
\usepackage{comment}
\usepackage{xcolor}

\usepackage{paralist}
\usepackage{mathtools}
\usepackage{colonequals}
\usepackage{amsmath}
\usepackage{amssymb}
\usepackage{amsthm}
\usepackage{tikz}
\usepackage{tikz-cd}
\usepackage{enumitem}
\usepackage{mleftright}
 
\usepackage{graphicx}
\usepackage{verbatim} 
\usepackage{setspace}
\usepackage[foot]{amsaddr}

\theoremstyle{plain}
\newtheorem{thm}{Theorem}[section]
\newtheorem{theorem}[thm]{Theorem}

\newtheorem{corollary}[thm]{Corollary}
\newtheorem{prop}[thm]{Proposition}
\newtheorem{proposition}[thm]{Proposition}

\newtheorem{lemma}[thm]{Lemma}

\newtheorem*{thm*}{Theorem}
\newtheorem*{theorem*}{Theorem}
\newtheorem*{lem*}{Lemma}
\newtheorem*{lemma*}{Lemma}
\newtheorem*{prop*}{Proposition}
\newtheorem*{proposition*}{Proposition}
\newtheorem*{cor*}{Corollary}
\newtheorem*{corollary*}{Corollary}
\newtheorem*{conj*}{Conjecture}
\newtheorem*{conjecture*}{Conjecture}

\theoremstyle{definition}
\newtheorem*{definition*}{Definition}
\newtheorem{definition}[thm]{Definition}
\newtheorem*{defn*}{Definition}

\theoremstyle{remark}
\newtheorem*{rem*}{Remark}
\newtheorem{rem}[thm]{Remark}
\newtheorem{remark}[thm]{Remark}
\newtheorem*{remark*}{Remark}
\newtheorem*{example*}{Example}
\newtheorem{example}[thm]{Example}

\newcommand{\acknowledgement}{\subsection*{Acknowledgements}}

\newcommand{\R}{\mathbb{R}}
\newcommand{\Z}{\mathbb{Z}}

\newcommand{\Ind}{\mathrm{Ind}}

\renewcommand{\colon}{\nobreak\mskip1mu\mathpunct{}\nonscript
  \mkern-\thinmuskip{:}\mskip3muplus1mu\relax}

\let\originalleft\left
\let\originalright\right
\renewcommand{\left}{\mathopen{}\mathclose\bgroup\originalleft}
\renewcommand{\right}{\aftergroup\egroup\originalright}

\title[Vietoris-Rips homology for semi-uniform spaces]{Vietoris-Rips homology
theory for semi-uniform spaces}

\date{}
\author{Antonio Rieser}
\address{Centro de Investigaci\'on en Matem\'aticas, Calle Jalisco S/N, Colonia Valenciana, C.P. 36023 Guanajuato, GTO, Mexico}
\thanks{Research supported in part by Cat\'{e}dras CONACYT / 1076 and Grant \#N62909-19-1-2134
from the US Office of Naval Research Global and the Southern Office of Aerospace Research and Development of the US Air Force Office of Scientific Research}

\begin{document}

\begin{abstract}
The Vietoris-Rips complex was first introduced in 1927 by Leopold Vietoris
\cite{Vietoris_1927}, and while it is now widely
used in both topological data analysis and the theory of hyperbolic groups,
many of the fundamental properties of the homology of the Vietoris-Rips complex have remained elusive. In this article, we
define the Vietoris-Rips homology for semi-uniform spaces, which generalizes the classical
theory for graphs and metric spaces and provides a natural, general setting for
the construction. We then prove a version of the Eilenberg-Steenrod axioms in
this setting, giving a natural definition of homotopy for semi-uniform
spaces in the process.\end{abstract}

\maketitle

\section{Introduction}

The Vietoris-Rips complex and its homology and cohomology have become
indispensible tools in the applications of algebraic topology,
particularly to areas with a combinatorial flavor. In topological
data analysis, it has become the standard
homology theory used to analyze data with persistent homology
\cites{Bauer_2019, Carlsson_Zomorodian_2005, Zomorodian_2010}, 
and before its emergence as a  
computational tool, it was already widely used in
the study of hyperbolic groups \cites{Ghys_de-la-Harpe_1990, Gromov_1987}.
Despite its popularity, however, little is known about even the most basic
properties of the Vietoris-Rips homology, except in the limiting case when the
defining parameter approaches zero, i.e. the so-called `metric cohomology'
studied by Hausmann in \cite{Hausmann_1995}. Indeed, the standard properties,
i.e. the Eilenberg-Steenrod axioms, fail for the classical Vietoris-Rips
homology with a fixed scale parameter $r > 0$, even in the form studied in
\cite{Hausmann_1995} for metric spaces, and to date, no alternative has been
proposed. 

In this article, we construct a version of the Vietoris-Rips homology by first
encoding the scale parameter into the structure of the space, and then
constructing a Vietoris-Rips homology theory for such space-scale pairs. Once
this is accomplished, we demonstrate that a version of the
Eilenberg-Steenrod axioms does, in fact, hold. The encoding of the scale parameter into the space is
achieved through the choice of a semi-uniform structure given by the parameter.
Indeed, we will develop the Vietoris-Rips theory for general semi-uniform spaces, and we do
not actually require the metric for the constructions. As a result, the theory
also applies to general topological space, graphs, and other cases of interest as well. In
the last section, we show that, on graphs and the $\leq$-version of the Vietoris-Rips
complex on finite metric spaces, our construction reduces to the standard
Vietoris-Rips theory, and the results in Section
\ref{sec:Eilenberg-Steenrod} are
new even for the classical Vietoris-Rips homology in these cases. We
additionally show that, as in the classical case, the homology of a Riemannian
manifold may be recovered by the (semi-uniform) Vietoris-Rips homology of a
sufficiently Gromov-Hausdorff close finite
metric space.

We remark, too, that studying the Vietoris-Rips theory from this point of view clarifies
that there is no essential difference between such `topological' constructions
for graphs and those for topological spaces. Indeed, these
results give an illustrative example of how considering homotopy and homology
in categories which strictly contain the categories of topological spaces, graphs, and
simplicial sets unifies many situations of interest, providing significant conceptual and
technical advantages over developing the relevant theory in each area separately.

\section{Semi-uniform spaces}

In this section, we give a brief introduction to semi-uniform spaces, which
will be our primary objects of study. Further details may be found in
\cite{Cech_1966}.

\subsection{Semi-uniform spaces}

Semi-uniform spaces and uniformly continuous functions between them generalize
the essential aspects of uniform continuity in metric spaces to a wider
collection of spaces and maps. In order to define semi-uniform structures, we will first need the following preliminary definitions.

\begin{definition}
	For a set $U \subset X \times X$, we define $U^{-1} \coloneqq \{(y,x)
	\mid (x,y) \in U\}$, and for $A\subset X$ and $U \subset X \times X$,
        we define \[U[A]
	\coloneqq \{y \in X \mid (a,y) \in U \subset X \times X \text{ for some
    } a \in A\}.\] If $A = \{x\}$, we write $U[x]$ for $U[\{x\}]$, and if $\mathcal{F}$ is a collection of subsets of $X \times X$ (such as a filter), then we define the collection of sets $[\mathcal{F}][A] \coloneqq \{U[A] \mid U \in \mathcal{F}\}$.
\end{definition}

We also recall the definition of a filter.

\begin{definition}
    \label{def:Filter}
    Let $X$ be a set, and let $\mathcal{F}$ be a non-empty collection of subsets of $A$. We say that $\mathcal{F}$ is a \emph{filter} on $X$ if
    \begin{enumerate}
	\item $\emptyset \notin \mathcal{F}$
	\item \label{def:Filter 2} $A, B \in \mathcal{F} \implies A \cap B \in \mathcal{F}$
	\item \label{def:Filter 3} $A \in \mathcal{F}, A \subset B \subset X \implies B \in \mathcal{F}$
    \end{enumerate}
\end{definition}

A collection $\mathcal{F}$ satisfying the above definition is sometimes called
a \emph{proper filter} in the literature, with \emph{filter} being reserved for
collections $\mathcal{F}$ which only satisfy items \ref{def:Filter 2} and
\ref{def:Filter 3} of Definition \ref{def:Filter}.
However, since none of the filters we will encounter contain the empty
set, we will not make this distinction in this article. 

We now proceed to the definition of a semi-uniform space.

\begin{definition}
    \label{def:Semi-uniform}
	Let $X$ be a set. A \emph{semi-uniform structure} on $X$ is a filter $\mathcal{F}$ on $X \times X$ such that 
\begin{enumerate}
\item \label{item:Diagonal}Each set in $\mathcal{F}$ contains the diagonal $\Delta \subset X \times X$
\item \label{item:Inverse set}If $U\in \mathcal{F}$, then $U^{-1}$ contains an element of $\mathcal{F}$
\end{enumerate} 
A \emph{semi-uniform space} $(X,\mathcal{F})$ is a pair consisting of a set $X$ and a semi-uniform structure $\mathcal{F}$ on $X$.
\end{definition}

\begin{remark}
    Note that, since $\mathcal{F}$ is a filter, then item \ref{item:Inverse
    set} in Definition \ref{def:Semi-uniform} above is equivalent to the
    condition
    \begin{enumerate}
        \item[(2')] $U^{-1} \in \mathcal{F}$.
    \end{enumerate}
\end{remark}

\begin{definition}
We say that a function $f:(X,\mathcal{F}_X) \to (Y,\mathcal{F}_Y)$ is \emph{uniformly continuous} iff, for each $V \in \mathcal{F}_Y$, there exists a $U \in \mathcal{F}_X$ such that $(x,y) \in U$ implies $(f(x),f(y)) \in V$.
\end{definition}

\begin{definition} Let $(X,\mathcal{F})$ be a semi-uniform space. If
    $\mathcal{G}$ is a filter base for $\mathcal{F}$, then we say that
    $\mathcal{G}$ is a \emph{base for the semi-uniformity $\mathcal{F}$}.
    Similarly, if $\mathcal{G}$ is a  subbase for the filter $\mathcal{F}$ then
    we say that $\mathcal{G}$ is a \emph{subbase for the semiuniformity
    $\mathcal{F}$}. 
\end{definition}

Our principal examples of semi-uniform spaces will be built from
semi-pseudometrics, defined here.

\begin{definition} A \emph{semi-pseudometric $d:X \times X \to \R$ on $X$} is a map such that
\begin{enumerate}
\item $d(x,x) = 0$ for all $x$ in $X$.
\item $d(x,y) = d(y,x)$ for all $(x,y) \in X \times X$.
\end{enumerate}
\end{definition}
\begin{example}
    \label{ex:Semi-uniform examples}
	\begin{enumerate}[wide]
		\item \label{ex:Semi-uniformity from metric} If $d$ is a semi-pseudometric
			on $X$, then the collection of sets $A_r = \{(x,y) \mid d(x,y)
< r\}$, where $r > 0$, form a filter base on $X \times X$ which satisfies conditions
(\ref{item:Diagonal}) and (\ref{item:Inverse set}) in Definition
\ref{def:Semi-uniform}. By Theorem 23.A.4 in \cite{Cech_1966}, this implies
that it is a filter base for a semi-uniform structure $\mathcal{F}$ on $X$,
which we will call the \emph{semi-uniformity induced by $d$.}
\item \label{ex:Semi-uniformity induced by d}Suppose that $(X,d)$ is a
    semi-psueometric
	space, and let $q\geq 0$. Then the function $d_q:X
	\times X \to \R$ defined by  
\begin{equation*}
d_q(x,y) = \begin{cases}
   d(x,y)-q & d(x,y) > q\\
   0 & d(x,y) \leq q
   \end{cases}
\end{equation*}
is also a semi-pseudometric on $X$. Denote by $\mathcal{F}_q$ the semiuniformity
induced by $d_q$ constructed as in
example \ref{ex:Semi-uniformity from metric} above. We will refer to
$\mathcal{F}_q$ as the \emph{semi-uniformity induced
by $d$ at scale $q$}.

\item \label{ex:Closed semi-uniformity induced by d}As above, let $(X,d)$ be a metric
	space, and suppose that $q > 0$. Define the function $d_{\leq q}:X
	\times X \to \R$ by 
\begin{equation*}
    d_{\leq q}(x,y) = \begin{cases}
	d(x,y) & d(x,y) > q\\
   0 & d(x,y) \leq q.
   \end{cases}
\end{equation*}
and note that $d_{\leq q}$ is a semi-psueodmetric on $X$. Denote by $\mathcal{F}_{\leq
q}$ the semi-uniformity
induced by $d_{\leq q}$ constructed as in
Example \ref{ex:Semi-uniformity from metric} above. We call $\mathcal{F}_{\leq q}$ is the \emph{closed semi-uniformity induced
by $d$ at scale $q$}. Note that the main difference between this structure and
the structure $\mathcal{F}_q$ in Example \ref{ex:Semi-uniformity induced by d} is that
\begin{equation*}
	\bar{A}_q \coloneqq \{ (x,y) \mid d(x,y) \leq q\} \in \mathcal{F}_{\leq q},
\end{equation*}
whereas $\bar{A}_q \notin \mathcal{F}_q$. 
\end{enumerate}

\end{example}

In the next proposition, we give a particularly nice characterization of uniformly continuous functions
$f:(X,\mathcal{F}_p) \to (Y,\mathcal{F}_q)$ between semi-uniform spaces
constructed as in Example
\ref{ex:Semi-uniform examples}(\ref{ex:Semi-uniformity induced by d}) above. We
begin with the following definition.
 
\begin{definition}
	Let $(X,d_X)$ and $(Y,d_Y)$ be semi-pseudometric spaces. We say that a function
	$f:X \to Y$ is $(p,q)$-continuous iff $\forall \epsilon>0$, $\exists
	\delta>0$ such that, for every $x,y \in X$ we have
	\begin{equation*}
		d_X(x,y) < p + \delta \implies d_Y(f(x),f(y)) < q + \epsilon.
	\end{equation*}
\end{definition}

\begin{proposition}
	Let $(X,d_X)$ and $(Y,d_Y)$ be semi-pseudometric spaces, and let $p,q \geq 0$. A
	function $f:X \to Y$ is $(p,q)$-continuous iff it is uniformly
	continuous as a function $f:(X,\mathcal{F}_p) \to (Y,\mathcal{F}_q)$
        between semi-uniform spaces.
\end{proposition}

\begin{proof}
	First, suppose $f:X \to Y$ is $(p,q)$-continuous, and let $V \in
	\mathcal{F}_q$.
	Then, by definition of $\mathcal{F}_q$, $\exists r>q$ such that
	\begin{equation}
		V_r \coloneqq \{(y,y') \in Y\times Y \mid d_Y(y,y') < r)\}
		\subset V.
	\end{equation} 
	Choose $0<\epsilon < r-q$. Since $f$ is $(p,q)$-continuous, there exists a
	$\delta >0$ such that $d_X(x,y) < p+\delta \implies d_Y(f(x),f(y)) < q+
	\epsilon$. Therefore,  $f(U_{p+\delta}) \subset V_r \subset V$, where
	\begin{equation*}
		U_{p+\delta} \coloneqq \{ (x,x') \in X \times X \mid d_X(x,x') <
		p + \delta \}.
	\end{equation*} 
	However, $U_{p+\delta} \in \mathcal{F}_p$ by definition of
	$\mathcal{F}_p$. Since $V \in
	\mathcal{F}_q$ was arbitrary, it follows that $f$ is uniformly continuous.

	Now suppose that $f$ is uniformly continuous as a map between the
	semi-uniform spaces
	$(X,\mathcal{F}_p) \to (Y,\mathcal{F}_q)$, and let $\epsilon >0$.
	Define
	\begin{equation*}
		V_\epsilon \coloneqq \{ (y,y') \in Y \times Y \mid d_Y(y,y') < q +
		\epsilon\}.
	\end{equation*}
	Since $f$ is uniformly continuous, and each $V_\epsilon \in
	\mathcal{F}_q$, then there exists a $U \in \mathcal{F}_p$
	such that $f(U) \subset V_\epsilon$. However, $\mathcal{F}_p$ is generated by
	sets of the form
	\begin{equation*}
		U_\delta \coloneqq \{ (x,x') \in X \times X \mid d_X(x,x') < p +
		\delta \},
	\end{equation*}
	where $\delta>0$. Therefore, there exists a $\delta > 0$ such that
	$U_\delta \subset U$, so we have $f(U_\delta) \subset f(U) \subset
	f(V_\epsilon)$, and it follows that $f$ is $(p,q)$-continuous.
\end{proof}

\subsection{Semi-uniform spaces from a \v Cech closure space}

In this section, we examine a special class of semi-uniform spaces which are constructed
from the interior covers of a \v Cech closure space $(X,c_X)$, and, in the case of a
topological space, from its open covers.
We begin
by giving a short introduction to \v Cech closure spaces. For more details on
\v Cech closure spaces, see
\cite{Cech_1966} and \cite{Rieser_2021}.

\begin{definition}
    Let $X$ be a set, and let $c:\mathcal{P}(X) \to \mathcal{P}(X)$ be a
    function on the power set of $X$ which satisfies
    \begin{enumerate}
        \item $c(\emptyset) = \emptyset$
        \item $A \subset c(A)$ for all $A \subset X$
        \item $c(A \cup B) = c(A) \cup c(B)$ for all $A,B \subset X$
    \end{enumerate}
    The function $c$ is called a \emph{\v Cech closure operator} (or \emph{closure
    operator}) on $X$, and the pair
    $(X,c)$ is called a \emph{\v Cech closure space} (or \emph{closure space}).
\end{definition}

\begin{example}  
    \begin{enumerate}[wide,labelwidth=!, labelindent=0pt]        
        \item Let $X$
            be a topological space with $A\subset X$, and denote by $\bar{A}$ the
            (topological) closure of $A$. Then $c(A) =
            \bar{A}$ is a \v Cech closure operator. Note that, in this case,
            $c^2(A) = c(A)$. For such closure operators (i.e. with $c^2=c$), it
            is known that the collection
            \begin{equation*}
                \mathcal{O} \coloneqq \{ X \backslash c(A) \mid A \subset X \}
            \end{equation*}
            forms the open sets of a topology on $X$ \cite{Cech_1966}. We will refer to closure operators with the
            property $c^2 = c$ as \emph{topological} or \emph{Kuratowski
            closure operators}, and we will denote by $\tau$ the topological closure structure of the
            topological space $(X,\mathcal{O})$. If the
            topology on $X$ is understood, then we
            write $(X,\tau)$ as the closure space with the topological closure
            structure $\tau$.
        \item Let $(X,d)$ be a metric space, and $r\geq 0$ a non-negative real number.
            For any $A \subset X$, define
            \begin{equation*}
                c_r(A) \coloneqq \{ x \in X \mid d(x,A) \leq r\}.
            \end{equation*}
            Then $c_r$ is a closure operator on $X$. Note that $c_0$ is the
            topological closure structure for the topology induced by the
            metric.
        \item Let $G = (V,E)$ be a graph. Then we define a closure structure
            $c_G$ on
            the set of vertices $V$ by 
            \begin{equation*}
                c_G(A) \coloneqq A \cup \{ w \in V \mid \exists v \in A : (v,w)
                \in E\} \text{ for all } A \subset V.
            \end{equation*}
            With this closure structure, $(V,c_G)$ is a closure space. 

            Note
            that $c_G$ has the property that 
            \begin{equation}
                \label{eq:Graph closure structure property}
                c_G(A) = \cup_{x \in A} c_G(x).
            \end{equation}
            If, conversely, we begin with a closure space $(X,c)$ such that $c$
            satisfies the property in Equation
            \ref{eq:Graph closure structure property}, then we may
            define a (directed) graph $G_X = (X,E_c)$ whose vertices are the points of $X$
            and whose edges are given by
            \begin{equation*}
                E_c \coloneqq \{ (v,w) \in X \times X \mid w \in c(v) \text{
                and } v \neq w\}.
            \end{equation*}
    \end{enumerate}
\end{example}

\begin{definition}
    Let $(X,c)$ be a closure space, and suppose $A\subset X$. We define the
    \emph{interior of A}, $i(A)$, by 
    \begin{equation*}
        i(A) = X - c(X - A).
    \end{equation*}
    Furthermore, we say that $U \subset X$ is a \emph{neighborhood of $A$} iff
    $A \subset i(U)$.
    A collection $\mathcal{U}$ of subsets of $X$ is said to be an \emph{interior cover} of
    $X$ iff \[X = \bigcup_{U \in \mathcal{U}} i(U).\]
\end{definition}

\begin{proposition}
    \label{prop:Interior of subsets}
    Let $(X,c)$ be a \v Cech closure space, and suppose that $B \subset A
    \subset X$. Then $i(B) \subset i(A)$.
\end{proposition}
\begin{proof}
    Since $X-A \subset X-B$, the axioms of the closure structure give $c(X-A) \subset c(X-B)$, and
    therefore \[i(B) = X - c(X-B) \subset X-c(X-A) = i(A),\] as desired.
\end{proof}

\begin{proposition}
    \label{prop:Intersection of neighborhoods}
    Let $(X,c)$ be a \v Cech closure space. For any two neighborhoods
$U$ and $V$ of a point $x \in X$, the intersection $U \cap V$ is also a
neighborhood of $x$. 
\end{proposition}
\begin{proof}
    Since $x \in i(U) \cap i(V)$, we have
    \begin{align*}
        x \in (X - c(X-U)) \cap (X - c(X-V))&= X - (c(X-U) \cup c(X-V))\\
        &= X - c((X-U) \cup (X-V))\\
        &= X - c(X - (U \cap V)),
    \end{align*}
    by the properties of $c$ and de Morgan's laws. Therefore, $x \in i(U\cap
    V)$.
\end{proof}

\begin{corollary}
    \label{cor:Intersection of interior covers}
    Suppose that $\mathcal{U}$ and $\mathcal{V}$ are interior covers of a \v
    Cech closure space $(X,c)$. Then the collection
    \begin{equation*}
        \mathcal{U} \cap \mathcal{V} \coloneqq \{ U \cap V \mid U \in
        \mathcal{U}, V \in \mathcal{V}\}
    \end{equation*}
    is an interior cover of $(X,c)$.
\end{corollary}
\begin{proof}
    Since $\mathcal{U}$ and $\mathcal{V}$ are interior covers of $(X,c)$, for
    every $x \in x$, there
exist $U_x \in \mathcal{U}$ and $V_x \in \mathcal{V}$ which are both
neighborhoods of $x$. By Proposition \ref{prop:Intersection of neighborhoods},
$U_x \cap V_x$ is a neighborhood of $x$ as well. Since $U_x \cap V_x \in
\mathcal{U} \cap \mathcal{V}$ and $x \in X$ is arbitrary, the conclusion
follows.
\end{proof}

We now begin the construction of two semi-uniform structures associated to a
\v Cech closure space. We first use an interior cover to define a pair of
relations which we will use to build the semi-uniform structures.

\begin{definition}
    \label{def:Relations from closure}
    Let $(X,c)$ be a \v Cech closure space, let $\mathcal{I}_c$ denote
    the collection of its interior covers of $(X,c)$, and let $\mathcal{U} \in
    \mathcal{I}_c$ be an
    interior cover of $(X,c)$. We define the relation $V_{\mathcal{U}}$, the \emph{Vietoris relation of the cover $\mathcal{U}$},
    by 
    \begin{equation*}
        V_\mathcal{U} = \{ (x,y) \mid \exists U \in \mathcal{U} \colon x,y \in
        U\}.
    \end{equation*}
    We additionally define the relation $R_\mathcal{U}$, the \emph{interior-inclusion (I-I) relation of the
    cover $\mathcal{U}$}, by
    \begin{equation*}
        R_{\mathcal{U}} = \{ (x,y) \mid \exists U \in \mathcal{U} \colon (x \in
        i(U) \text{ and } y \in U)\text{ or }(y \in i(U) \text{ and } x \in U)\}.
    \end{equation*}    
    We denote by $\mathcal{V}_c$ and
    $\mathcal{R}_c$ the filters on $X \times X$
generated by the collections $\{V_{\mathcal{U}}\}_{\mathcal{U} \in
\mathcal{I}_c}$ and $\{R_{\mathcal{U}}\}_{\mathcal{U} \in
\mathcal{I}_c}$, respectively.

\end{definition}

We will now show that $\mathcal{V}_c$ and $\mathcal{R}_c$ define semi-uniform structures
on $X$. We first
require the following lemma.

\begin{lemma}
    \label{lem:Refinement of covers}
    Let $\mathcal{U}$ and $\mathcal{U}'$ be interior covers of the \v
    Cech closure space $(X,c)$. If $\mathcal{U} < \mathcal{U}'$ (i.e. if
    $\mathcal{U}'$ refines $\mathcal{U}$), then $V_{\mathcal{U}'} \subset
    V_{\mathcal{U}}$ and $R_{\mathcal{U}'} \subset R_{\mathcal{U}}$.
\end{lemma}
\begin{proof}
    We first show that $R_{\mathcal{U}'} \subset
    R_{\mathcal{U}}$. Suppose $(x,y) \in R_{\mathcal{U}'}$. Then there exists a
    $U' \in
    \mathcal{U}'$ such that, without loss of generality, $x \in i(U')$ and $y
    \in U'$. 
    However, since $\mathcal{U}'$ refines $\mathcal{U}$, there exists a set $U
    \in \mathcal{U}$ 
    with $U' \subset U$. By Proposition \ref{prop:Interior of subsets}, $i(U')
    \subset i(U)$ as
    well. It follows that $x \in i(U)$ and $y \in U$, and therefore $(x,y) \in
    R_\mathcal{U}$, as desired.

    Now suppose that $(x,y) \in V_{\mathcal{U}'}$. Then there exists a $U'
    \in \mathcal{U}'$ with $x,y \in U'$. As above, since $\mathcal{U}'$ refines
    $\mathcal{U}$, there exists a $U \in \mathcal{U}$ with $U' \subset U$, and
    therefore $x,y \in U$, which gives $(x,y) \in V_{\mathcal{U}}$, proving the
    result.
\end{proof}

\begin{proposition}
    Let $(X,c)$ be a closure space. Then $(X,\mathcal{V}_c)$ and
    $(X,\mathcal{R}_c)$ are semi-uniform spaces.
\end{proposition}
\begin{proof}
    By 23.A.4 in \cite{Cech_1966}, for a filter $\mathcal{F}$ to define a
    semi-uniform structure on $X$, we must have that
    \begin{enumerate}
        \item \label{item:DiagonalCondition} Every element $U \in \mathcal{F}$ contains the diagonal in $X
            \times X$, and
        \item \label{item:Inverse} For any element $U \in \mathcal{F}$, 
            the inverse $U^{-1}$ contains an element of
            $\mathcal{F}$.
    \end{enumerate}
    We first show that $\mathcal{V}_c$ satisfies these conditions. Let
    $V_\mathcal{U} \in \mathcal{V}_c$. Then, by construction, for every $x\in X$
    exists $U \in \mathcal{U}$ such that $x \in U$. It follows that $(x,x) \in
    V_\mathcal{U}$, proving Item (\ref{item:DiagonalCondition}). Furthermore, $(x,y) \in
    V_\mathcal{U}$ iff there exists a $U \in \mathcal{U}$ with $x,y \in U$,
    which implies that $(y,x) \in V_\mathcal{U}$ as well, proving Item 
    (\ref{item:Inverse}). 

    Now let $R_\mathcal{U} \in \mathcal{R}_c$. Then, as above, for any $x \in X$,
    there is a $U \in \mathcal{U}$ with $x \in i(U)$, and therefore $(x,x) \in
    R_\mathcal{U}$, fulfilling Item (\ref{item:DiagonalCondition}) above. Also, by construction, $R_\mathcal{U}^{-1} =
    R_\mathcal{U}$, and so Item (\ref{item:Inverse}) is satisfied as well,
    and the proof is complete.
\end{proof}

\begin{definition}
    \label{def:Topological semi-uniform structure}
    Given a \v Cech closure space $(X,c)$, we say that the pair $(X,\mathcal{V}_c)$ is
    the \emph{Vietoris semi-uniform space induced by $c$}, and that the pair
    $(X,\mathcal{R}_c)$ is the \emph{interior-inclusion} (or \emph{I-I}) \emph{semi-uniform space induced by
    $c$}.
\end{definition}

\begin{proposition}
Let $(X,c)$ be a topological closure space (i.e. where $c^2 = c$). Then $\mathcal{V}_c = \mathcal{R}_c$. 
\end{proposition}

\begin{proof} Denote by $\mathcal{O}$ the set of open covers of $X$, and let
    $\mathcal{U} \in \mathcal{O}$. Since, for every $U
    \in \mathcal{U}$, $i(U) = U$, we have that $V_\mathcal{U} =
    R_\mathcal{U}$. Since any interior cover
$\mathcal{U}'$ of $X$ is refined by an open cover $\mathcal{U} \in
\mathcal{O}$, it follows from Lemma \ref{lem:Refinement of
covers} that the sets $\{V_\mathcal{U}\}_{\mathcal{U} \in \mathcal{O}}$ and
$\{R_\mathcal{U}\}_{\mathcal{U} \in \mathcal{O}}$ generate the filters
$\mathcal{V}_c$ and $\mathcal{R}_c$, respectively. Since
$\{V_\mathcal{U}\}_{\mathcal{U} \in \mathcal{O}} =
\{R_\mathcal{U}\}_{\mathcal{U} \in \mathcal{O}},$
the conclusion follows.
\end{proof}

\begin{definition} Let $(X,c)$ be a topological closure space as above. We call the
    semi-uniform structure $\mathcal{V}_c = \mathcal{R}_c$ the
    \emph{topological semi-uniform structure induced by $c$}, which we will denote
    by $\mathcal{T}_c$. We call the
    semi-uniform space $(X,\mathcal{T}_c)$ the \emph{topological
    semi-uniform space induced by $c$}.\end{definition}

\section{The Vietoris-Rips Complex and Vietoris-Rips Homology}
\label{sec:VR homology construction}

Let $(X,\mathcal{F})$ be a semi-uniform space. For every element $U \in
\mathcal{F}$, we now define a simplicial
complex from the pair $(X,U)$, which we denote $\Sigma^X_U$.  We proceed as follows. First, we construct the graph $G(X,U) = (V,E)$ by:
\begin{align*}
V &= \{x \mid x \in X \}\\
E &= \{ (x_0,x_1) \mid (x_0,x_1) \in U\}
\end{align*}

	Note that, a priori, this gives a directed graph. If, for every
	$(x_0,x_1)$, the edge $(x_1,x_0)$ also exists, then we replace the
	directed graph with an undirected graph, which we still denote
	$G(X,U)$. If $G(X,U)$ is undirected, we let $\Sigma^X_U$ be the clique
	complex of $G(X,U)$. If $G(X,U)$ is directed, then we let
	$\Sigma^{X}_{U}$ be the directed clique complex of $G(X,U)$. That is,
	in the directed case, we say that a finite set $\{x_0,\dots,x_q\}$ is a
	$q$-simplex in $\Sigma^{X}_{U}$ iff there exists a permutation of the
	indices $\gamma:\{0,\dots,q\} \to \{0,\dots,q\}$ such that
	$(x_{\gamma(i)},x_{\gamma(j)}) \in E$ if $\gamma(i) < \gamma(j)$. Note
	that we do not assume a global ordering of the vertices, and so all
	possible orderings of the vertices $\{x_0,\dots,x_n\}$ are considered
	when determining whether the collection forms a simplex in the directed
	flag complex. If there exists such a permutation, then for any subset
	of $\{x_{\gamma(0)},\dots,x_{\gamma(q)}\}$ the condition
	$(x_{\gamma(i)},x_{\gamma(j)}) \in E$ if $\gamma(i) < \gamma(j)$ is
	also fulfilled, and so $\Sigma^{X}_{U}$ is indeed a simplicial complex.

\begin{definition}
	We call the simplicial complex $\Sigma^X_U$ the \emph{Vietoris-Rips complex of the relation $U \subset X \times X$}. 
\end{definition}

\subsection{Vietoris-Rips homology and cohomology}

Let $(X,\mathcal{F})$ be a semi-uniform space. For $A\subset X$ and $U \in
\mathcal{F}$, define
\begin{equation*}
	U_A \coloneqq U \cap (A \times A),
\end{equation*}
and let 
$\mathcal{F}_A$ be the collectino of sets
\begin{equation*}
	\mathcal{F}_A \coloneqq \{ U_A \mid U \in \mathcal{F}, U_A \neq
	\emptyset \},
\end{equation*}
which we call the \emph{relativation of $\mathcal{F}$ by $A$}.
\begin{proposition}
	Let $(X,\mathcal{F})$ be a semi-uniform space. For any subset $A \subset X$, the relativization $\mathcal{F}_A$ forms a semi-uniform structure on $A$.
\end{proposition}
\begin{proof}
	Immediate from the definitions.
\end{proof}

\begin{definition}
	Let $A \subset X$. We write $(X,A,\mathcal{F})$ to denote the pair of semi-uniform spaces
	$((X,\mathcal{F}),(A,\mathcal{F}_A))$, which we call a \emph{semi-uniform pair}. 
\end{definition}
The relation between the Vietoris-Rips
complexes $\Sigma^A_{U_A}$ and $\Sigma^X_{U}$ is given by the following
proposition. 

\begin{prop}
	Let $(X,\mathcal{F})$ be a semi-uniform space. Suppose that $A \subset
	X$, and $U \in \mathcal{F}$. Then $\Sigma^A_{U_A}$ is a
	subcomplex of $\Sigma^{X}_{U}$.
\end{prop}
\begin{proof}
	Let $\sigma \coloneqq (x_0,\dots,x_n) \in \Sigma^A_{U_A}$, and without loss of generality suppose
that $(x_i,x_j) \in U_A$ for all $i<j$. Then, in particular, $(x_i,x_j) \in U$
for all $i<j$, and therefore $\sigma \in \Sigma^{X}_{U}$, as desired.
\end{proof}

\begin{definition}
	Let $(X,\mathcal{F})$ be a semi-uniform space, $A\subset X$, $U\in
	\mathcal{F}$, and let $\Sigma^X_U$ and
	$\Sigma^A_{U_A}$ be the respective Vietoris-Rips complexes. We denote
	the homology and cohomology of the simplicial pair
	$(\Sigma^X_U,\Sigma^A_{U_A})$ with coefficients in
	the abelian group $G$ by
	\begin{align*}
		H^{VR}_*(X,A\colon U; G) &\coloneqq
		H_*\left(\Sigma^X_U,\Sigma^A_{U_A}; G\right)\\
		H^*_{VR}(X,A\colon U; G) &\coloneqq H^*(\Sigma^X_U,\Sigma^A_{U_A}; G).
	\end{align*}
We call these groups the \emph{Vietoris-Rips homology} and the \emph{Vietoris-Rips
cohomology of $(X,A,U)$}, repsectively. When $A = \emptyset$, we will simply
write $H_*^{VR}(X\colon U;G)$ and $H_{VR}^*(X\colon U;G)$. Additionally, when
$G = \Z$, we write $H^{VR}_*(X,A\colon U)$ and $H_{VR}^*(X,A\colon U)$ for
$H_*^{VR}(X,A\colon U; \Z)$ and $H^*_{VR}(X,A\colon U;\Z)$, respectively.

\end{definition}

\begin{lemma}
\label{lem:Homomorphism from relations}
Let $(X,A,\mathcal{F})$ be a semi-uniform pair, and let
$U, U' \in \mathcal{F}$ with $U \subset U'$. Then there is a simplicial embedding 
\[\iota_{U'U}: (\Sigma^X_{U},\Sigma^A_{U_A}) \to
(\Sigma^X_{U'},\Sigma^A_{U'_A})\] which induces homomorphisms on the simplicial homology and cohomology groups
\begin{align*}
\iota_*^{U'U}\colon H_*^{VR}(X,A\colon U;G) \to H_*^{VR}(X,A\colon U';G)\\
\iota^*_{UU'}\colon H^*_{VR}(X,A\colon U';G) \to H^*_{VR}(X,A\colon U;G).
\end{align*}
\end{lemma}

\begin{proof} Define $\iota_{U'U}$ by
    \begin{equation*}
        \iota_{U'U}(x_0,\dots,x_n) = (x_0,\dots,x_n)
    \end{equation*}
    Since $U \subset U'$, any simplex $\sigma \in \Sigma^X_U$ is also in $\Sigma^X_{U'}$, and it
    follows that $\iota_{U'U}$ is a simplicial embedding. Restriction of
    $\iota_{U'U}$ to $\Sigma^A_{U_A}$ gives the corresponding injection, and the
    homomorphisms on simplicial homology and cohomology are produced by the
    functoriality of homology and cohomology of simplicial complexes.
\end{proof}

Define a partial order $<$ on $\mathcal{F}$ by: $U <
U'$ iff $U \subset U'$. Since $\mathcal{F}$ is closed under finite intersections by
definition, $(\mathcal{F},<)$ is a directed set,
and Lemma \ref{lem:Homomorphism from relations} makes the collections of homology and cohomology
groups $\{H_*^{VR}(X,A\colon U;G)\}_{U\in \mathcal{F}}$ and $H^*_{VR}(X,A\colon U; G)$ into inverse and directed
systems, respectively. This observation leads to the following definitions.

\begin{definition}
	We define the \emph{Vietoris-Rips homology of $(X,A,\mathcal{F})$} with
	coefficients in $G$ by \[H^{VR}_*(X,A,\mathcal{F};G) = \varprojlim H^{VR}_*(X,A\colon U;G). \] 
	The \emph{Vietoris-Rips cohomology of $(X,A,\mathcal{F})$} is given by
	\[H_{VR}^*(X,A,\mathcal{F};G) = \varinjlim H_{VR}^*(X,A\colon U,\mathcal{U};G).  \]
\end{definition}

\begin{remark}
	 Note that, since $U$ and $U^{-1}$ are in $\mathcal{F}$, and since
 $\mathcal{F}$ is a (proper) filter, we have that $\emptyset \neq U \cap U^{-1}
 \in \mathcal{F}$, and so the subcollection of symmetric elements of
 $\mathcal{F}$ is cofinal in $\mathcal{F}$. This allows us to compute the
 Vietoris-Rips homology and cohomology of $(X,A,\mathcal{F})$ solely from the clique complexes of undirected graphs. 
 \end{remark}

 \begin{remark} Let $(X,c)$ be a \v Cech cloure space. We remark that, by results in \cite{Dowker_1952}, the
     Vietoris-Rips cohomology of $(X,\mathcal{V}_c)$ is isomorphic to the \v
     Cech cohomology of $(X,c)$ as defined in \cite{Palacios_Lic_2019} and
     \cite{Demaria_Garbaccio_1984}, i.e. the limit of the cohomologies of the nerve
     complexes of the interior covers, partially ordered by refinement. See
     \cite{Palacios_Lic_2019} for more details on the properties of the \v Cech cohomology of a closure space.
\end{remark}

\section{The Eilenberg-Steenrod axioms}
\label{sec:Eilenberg-Steenrod}
The Eilenberg-Steenrod axioms form the essential elements of a homology theory
on topological spaces. In this section, we give a version of the Eilenberg-Steenrod axioms for
semi-uniform spaces, and we show that they are satisfied by the Vietoris-Rips 
homology and cohomology defined in Section \ref{sec:VR homology construction}
above.

\subsection{Functoriality}

We first construct the homomorphisms on homology and cohomology induced by
uniformly continuous maps of the semi-uniform pairs \[f\colon (X,A,\mathcal{F}_X) \to
(Y,B,\mathcal{F}_Y),\] i.e. where the map $f\colon (X,\mathcal{F}_X) \to
(Y,\mathcal{F}_Y)$ is uniformly continuous and $f(A) \subset B$. We begin with the following lemma.

\begin{lemma}
\label{lem:Functoriality to simplicial complexes}
Let $f\colon (X,A,\mathcal{F}_X) \to (Y,B,\mathcal{F}_Y)$ be a uniformly continuous
map of semi-uniform pairs. Then, for every $V\in \mathcal{F}_Y$, there exists a $U \in
\mathcal{F}_X$ and an induced simplicial map of simplicial pairs
\[\tilde{f}_{VU}\colon \left(\Sigma^X_U,\Sigma^A_{U_A}\right) \to
\left(\Sigma^Y_V,\Sigma^B_{V_B}\right).\] 
\end{lemma}

\begin{proof}
We first note that, since $f$ is uniformly continuous, then, by definition, for every
$V \in \mathcal{F}_Y$, there exists a $U \in \mathcal{F}_X$ such that $f(U)
\subset V$. We construct the induced map $\tilde{f}:\Sigma^X_U \to
\Sigma^Y_V$. First, without loss of generality, let $\sigma=(x_0,\dots,x_n) \in
\Sigma^X_U$ be such that $(x_i,x_j) \in U$ for all $i < j$. Then define 
\begin{equation*}
    \tilde{f}(\sigma) \coloneqq (f(x_0),\dots,f(x_n)).\\
\end{equation*}
Since $f(U) \subset V$, it follows that $(f(x_i),f(x_j)) \in V$ for any $i <
j$, and therefore $\tilde{f}(\sigma) \in \Sigma^Y_V$, making $\tilde{f}$ a
simplicial map. Furthermore, since $f(A) \subset B$,
$\tilde{f}\left(\Sigma^A_{U_A}\right)
\subset \Sigma^B_{V_B}$, and therefore $\tilde{f}$ is a simplicial map of pairs, as
desired.
\end{proof}

\begin{lemma}	
\label{lemma:Commuting maps for functoriality}
Let $(X,\mathcal{F}_X)$ and $(Y,\mathcal{F}_Y)$ be semi-uniform spaces and
suppose the map $f:(X,\mathcal{F}_X) \to (Y,\mathcal{F}_Y)$ is uniformly
contiuous. Let $V, V' \in \mathcal{F}_Y$ with $V \subset V'$, and let $U,U'\in
\mathcal{F}_X$ be such that $U \subset U'$, $f(U) \subset V$, and $f(U') \subset V'$.
Then the induced maps on homology and cohomology make the diagrams
\begin{equation*}
\begin{tikzcd}
    H_*(X,A\colon U;G) \ar [r,"\tilde{f}_*^{VU}"] \ar [d,"\iota^*_{U'U}"] &
H^*(Y,B\colon V;G) \ar [d,"\iota^*_{V'V}"] \\
    H_*(X,A\colon U';G) \ar [r,"\tilde{f}_*^{V'U'}"] & H_*(Y,B\colon V';G)
\end{tikzcd}
\end{equation*} 
and
\begin{equation*}
\begin{tikzcd}
    H^{*}(Y,B \colon V';G) \ar [r,"\tilde{f}^*_{VU}"] \ar [d,"\iota^*_{VV'}"] 
 & H^{*}(X,A\colon U';G)     \ar [d,"\iota^*_{UU'}"] \\
    H^*(Y,B\colon V;G) \ar [r,"\tilde{f}^*_{V'U'}"] & H^*(X,A\colon U; G)
\end{tikzcd}
\end{equation*}
commute.
\end{lemma}
\begin{rem}
    Note that the existence of $U,U' \in \mathcal{F}_X$ is guaranteed by the
    uniform continuity of $f$.
\end{rem}
\begin{proof}
	We first consider the diagram of simplicial complexes
\begin{equation*}
\begin{tikzcd}
    \Sigma^X_U \ar [r, "\tilde{f}_{VU}"] \ar [d, "\iota_{U'U}"] &
\Sigma^Y_V \ar [d, "\iota_{V'V}"]\\
    \Sigma^X_{U'} \ar [r, "\tilde{f}_{V'U'}"] & \Sigma^Y_{V'}
\end{tikzcd}
\end{equation*}
where the $\iota_{**}$ are defined as in Lemma \ref{lem:Homomorphism from
relations} and the $\tilde{f}_{**}$ are the simplicial maps from Lemma
\ref{lem:Functoriality to simplicial complexes}. Let $\sigma = (x_0,\dots,x_n)
\in \Sigma^X_U$. By construction, $\iota_{U'U}(\sigma) = \sigma$, and
$\tilde{f}_{V'U'}(\sigma) = (\tilde{f}(x_0),\dots,\tilde{f}(x_n))$. Similarly, we 
have that $f_{VU}(\sigma) = (f(x_0),\dots,f(x_n))$, and
$\iota_{V'V}(f(x_0),\dots,f(x_n)) = (f(x_0),\dots, f(x_n))$. Since the 
 simplex $\sigma\in \Sigma^X_U$ was arbitrary, the diagram commutes. Applying the homology and 
 cohomology functors to this diagram, and the result follows.
\end{proof}

\begin{theorem}
\label{thm:Functoriality}
Let $(X,A,c)$ be a closure space pair, where $(X,c)$ is semi-uniformizable. The Vietoris-Rips homology $H_*^{VR}(X,A,c)$ and cohomology $H^*_{VR_k}(X,A,c)$ are functorial with respect to uniformly continuous maps.
\end{theorem}
\begin{proof}
The theorem follows from Lemmas \ref{lem:Functoriality to simplicial complexes}
and \ref{lemma:Commuting maps for functoriality}, taking direct limits.
\end{proof}

\subsection{Exactness}

\begin{theorem}
\label{thm:Exactness}
Let $(X,A,\mathcal{F})$ be semi-uniform pair, and $G$ an abelian group. Then there is a functorial long-exact sequence of pairs
\begin{equation*}
\cdots \to H^*_{VR}(X,A;G) \to H^*_{VR}(X;G) \to H^*_{VR}(A; G) \to
H^{*+1}_{VR}(X,A;G) \to \cdots
\end{equation*}
and a functorial sequence of order two
\begin{equation*}
\cdots \to H_*^{VR}(A;G) \to H_*^{VR}(X;G) \to H_*^{VR}(X,A;G) \to
H_{*-1}^{VR}(A;G) \to \cdots
\end{equation*}
\end{theorem}

\begin{proof}
    For each element $U\in \mathcal{F}$, one has the corresponding long exact
    sequences from simplicial homology and cohomology. Since exactness is
    preserved by direct limits, and exactness is of order two for inverse
    limits, it follows the sequence of pairs is exact for $H^*_{VR}(X,A;G)$ and
    of order two for $H_*^{VR}(X,A;G)$.
\end{proof}

\subsection{The Dimension Axiom}

\begin{theorem}
\label{thm:Dimension}
Let $(p,\mathcal{F}_p)$ denote the one point semi-uniform space. Then 
\[
    H_n^{VR}(p;G) \cong H^n_{VR}(p;G)\cong \begin{cases} 
G & n = 0\\
0 & n > 0,
\end{cases}
\]
\end{theorem}
\begin{proof}
Since the one-point set only has a single point, the dimension axiom follows from the corresponding result for simplicial complexes. 
\end{proof}

\subsection{The Excision axiom}

For the excision axiom, we first note that the simplicial homology and
cohomology groups of a combinatorial simplicial complex are isomorphic to the homology
and cohomology groups of its geometric realization. We first quote the following
result from \cite{Switzer_1975}, which we will use in the proof of the lemma which
follows.

\begin{lemma}[\cite{Switzer_1975}, Proposition 7.5]
\label{lem:Excision CW complex}
If the CW-complex $\Sigma$ is the union of subcomplexes $\Sigma_1$ and
$\Sigma_2$, then
the inclusion $i:(\Sigma_1, \Sigma_1 \cap \Sigma_2) \hookrightarrow (\Sigma,
\Sigma_2)$ induces an isomorphism
\[
    i_n:H_n(\Sigma_1, \Sigma_1 \cap \Sigma_2) \xrightarrow{\cong} H_n(\Sigma, \Sigma_2).
\]
for all $n \in \Z$.
\end{lemma}

\begin{remark}
    \label{rem:Excision duality remark}
    As mentioned on page 125 of \cite{Switzer_1975}, the dual of the above
    result also holds, by dualizing its proof. Therefore, we also have that $i$
    induces an isomorphism in cohomology
    \begin{equation*}
        i^n:H^n(\Sigma,\Sigma_2)\xrightarrow{\cong} H^n(\Sigma_1,\Sigma_1 \cap \Sigma_2) 
    \end{equation*}
    for every $n \in \Z$ as well.
\end{remark}
\begin{theorem}[Excision]
\label{thm:Excision}
Let $(X,\mathcal{F})$ be a semi-uniform space, and let $B \subset A \subset X$.
Suppose that there exists an element $W \in \mathcal{F}$
such that, for every $U \subset W$, $U \in \mathcal{F}$,
\[U[B] \subset A.\] 
Then
\begin{align*}
    H^*_{VR}(X-B,A-B:\mathcal{F};G) &\cong H^*_{VR}(X,A:\mathcal{F};G)\\
                                    \text{
                                    and }&\\
    H_*^{VR}(X-B,A-B:\mathcal{F};G) &\cong H^*_{VR}(X,A:\mathcal{F};G).
\end{align*}
\end{theorem}

\begin{proof}
    Let $A \subset B \subset X$, and let $U' < W, U' \in \mathcal{F}'$. Let $U
    = U' \cap
    (U')^{-1}$, and note that $U$ is symmetric, i.e. $U = U^{-1}$. We first show that
    \[\Sigma^X_U= \Sigma^A_{U_A} \cup \Sigma^{(X-B)}_{U_{(X-B)}}.\] 
    Let $\sigma = (x_0,\dots,x_n) \in \Sigma^X_U$ and suppose that $\sigma
    \notin \Sigma^{(X-B)}_{U_{(X-B)}}$. Then there must be at least one $x_i \in
    \{x_0,\dots,x_n\}$ with $x_i \in B\subset A$. However, since $U$ is symmetric, this
    implies that $x_j \in U[x_i] \subset U[B]\subset A$ for all $x_j \in \sigma$. It
    follows that $(x_j,x_k) \in U \cap (A \times A)$ for all $x_j,x_k \in
    \{x_0,\dots,x_n\}$, and therefore 
    \[\sigma \in
    \Sigma^A_{U_A}.\]

    Now suppose that $\sigma \notin \Sigma^A_{U_A}$. Then there exists $x_i \in
    \{x_0,\dots,x_n\}$ such that $x_i \in X - A$. Suppose there exists a $j$
    where $x_j \in B$. Then $x_i \in
    U[x_j] \subset U[B] \subset A$, and therefore $x_i \in  A \cap (X-A) =
    \emptyset$, a contradiction. Therefore, $x_j \in X-B$ for all
    $x_j \in \sigma$. Furthermore, we have $(x_j,x_k) \in U \cap ((X - B)
    \times (X-B))$
    for all $x_j,x_k \in \{x_0,\dots,x_n\}$, and therefore 
    \[ \sigma \in \Sigma^{(X-B)}_{U_{(X-B)}}. \]

    Taking $\Sigma_1 = X-B$ and $\Sigma_2 = A$ in Lemma \ref{lem:Excision CW
    complex} and Remark \ref{rem:Excision duality remark}, we have
    \begin{align*}
        H^*_{VR}(X-B,A-B:U;G) &\cong H^*_{VR}(X,A:U;G)\\ \text{ and }& \\
        H_*^{VR}(X-B,A-B:U;G) &\cong H^*_{VR}(X,A:U;G).
    \end{align*}
    Since the set of symmetric controlled sets $U$ is cofinal in $\mathcal{F}$, the result follows.
\end{proof} 

\subsection{The Homotopy axiom}
For the Homotopy axiom, we will follow the general line of proof in the
demonstration of the corresponding result for \v Cech homology and cohomology in
\cite{Eilenberg_Steenrod_1952}. Before we begin, we define several terms.

\begin{definition}
    Let $I \coloneqq [0,1]$, and let $\mathcal{F}_\tau$
    be the topological semi-uniform structure induced by the standard topology
    on $I$ as in Definition \ref{def:Topological semi-uniform structure}. Let $(X,\mathcal{F}_X)$ be
    a semi-uniform space. We write $(X \times I,\mathcal{F}_\Pi)$ for the
    product of the
    semi-uniform spaces $(X,\mathcal{F}_X)$ and $(I,\mathcal{F}_\tau)$\end{definition}

\begin{definition}
    Suppose that $(X,\mathcal{F}_X)$ and $(Y,\mathcal{F}_Y)$
    are semi-uniform spaces. 
We say that two uniformly continuous maps $f,g\colon (X,\mathcal{F}_X)
    \to (Y,\mathcal{F}_Y)$ are \emph{homotopic}, or $f \simeq g$, iff there exists a uniformly
    continuous map 
    \[ H:(X \times I,\mathcal{F}_\Pi) \to (Y,\mathcal{F}_Y)\] 
    such that 
    \[
        \begin{cases} H(x,0) = f(x) \\
                      H(x,1) = g(x)
        \end{cases}
    \]
    for all $x \in X$.

    We say that two semi-uniform spaces $(X,\mathcal{F}_X)$ and
    $(Y,\mathcal{F}_Y)$ are \emph{homotopy equivalent}, or $(X,\mathcal{F}_X)
    \simeq (Y,\mathcal{F}_Y)$, iff there exist uniformly continuous maps
    $f\colon (X,\mathcal{F}_X) \to (Y,\mathcal{F}_Y)$ and $g\colon (Y,\mathcal{F}_Y)
    \to (X,\mathcal{F}_X)$ such that $g\circ f \simeq 1_X$ and $f\circ g \simeq
    1_Y$. 
\end{definition}

We show that the topological semi-uniform structure on the interval $I=[0,1]$ has a more convenient
form, which will be useful in the following.

\begin{proposition}
    \label{prop:Metric and topological structures on compact spaces}
    Let $(X,d)$ be a compact metric space with topological closure operator
    $\tau$ (i.e. $\tau^2 = \tau$) induced by the metric. Then $\mathcal{F}_0 = \mathcal{F}_\tau$,
    where $\mathcal{F}_\tau$ is the topological semi-uniform structure induced
    by $\tau$, and $\mathcal{F}_0$ is the semi-uniform structure generated by
    the sets $\{(x,y) \mid d(x,y) < r\}$ where $r>0$ ranges over the positive
    real numbers, as in Example
    \ref{ex:Semi-uniform examples}(\ref{ex:Semi-uniformity from metric}).  
\end{proposition}

\begin{proof}
    We claim that $\mathcal{F}_\tau$ is generated by the collection of sets of the form $U_r
    \coloneqq \{(x,y) \mid
    d(x,y) < r \}$, where $r>0$ ranges over the positive real numbers. 
    We first show that every $U_r \in \mathcal{F}_\tau$. Fix $r>0$ and consider the open cover $\mathcal{B}_r \coloneqq
    \{B_x(r) \mid x \in X\}$, where
    $B_x(r)$ denotes the open ball centered at $x$ with
    radius $r$. Recall that, for the open cover $\mathcal{B}_r$, the Vietoris
    relation $V_{\mathcal{B}_r}$ is given by
    \[V_{\mathcal{B}_r} \coloneqq \{(y,y') \mid \exists x\in X \text{ such that
} y,y' \in B_x(r)\}.\] It follows that $U_r \subset V_{\mathcal{B}_r} \subset U_{2r}$. Since $V_{\mathcal{B}_r} \in
    \mathcal{F}_\tau$ and $\mathcal{F}_\tau$ is a filter, we have that $U_{2r}
    \in \mathcal{F}_\tau$, and since $r>0$ was arbitrary, $U_r \in
    \mathcal{F}_\tau$ for every $r>0$.

    Now suppose that $\mathcal{U}$ is an open cover of 
    $X$, and let $V_\mathcal{U}$ denote the corresponding Vietoris relation. 
    Since $X$ is compact by hypothesis, $\mathcal{U}$ has a Lebesgue number $\lambda$, and
    therefore the collection $\mathcal{B}_{\lambda}\coloneqq
    \{B_x(\lambda) \mid x \in X\}$ is an
    open cover of $X$ which refines $\mathcal{U}$. Since $\mathcal{U} <
    \mathcal{B}_{\lambda}$, from Proposition \ref{lem:Refinement of
    covers}, we have that $U_\lambda \subset V_{\mathcal{B}_\lambda} \subset V_{\mathcal{U}}$,
    proving the claim. 

    Since both filters $\mathcal{F}_\tau$ and $\mathcal{F}_0$ are generated by the same
    collection of sets, it follows that $\mathcal{F}_\tau = \mathcal{F}_0$.
\end{proof}

The remainder of this section will be devoted to proving the homotopy invariance of the Vietoris-Rips
homology and cohomology, i.e.

\begin{theorem}
\label{thm:Homotopy}
    Let $(X,\mathcal{F}_X)$ and $(Y,\mathcal{F}_Y)$ be semi-uniform spaces,
    and suppose $f_0,f_1:(X,\mathcal{F}_X) \to (Y,\mathcal{F}_Y)$ are
    uniformly continuous maps with $f_0 \simeq f_1$. Then
    \begin{align*}
        f_{0*} &= f_{1*}:H_*^{VR}(X,A,\mathcal{F}_X) \to H^{VR}_*(Y,B,\mathcal{F}_Y),
        \text{ and }\\
        f^*_0 &= f^*_1:H^*_{VR}(Y,B,\mathcal{F}_Y) \to H_{VR}^*(X,A,\mathcal{F}_X)
    \end{align*}	 
\end{theorem}

We first show that it will be sufficient to prove the following theorem.

\begin{theorem}
	\label{thm:Homotopy 2}
	Let $g_0,g_1: (X,\mathcal{F}_X) \to (X,\mathcal{F}_X) \times (I,\tau)$
        be defined by $g_0(x)=(x,0)$, $g_1(x) = (x,1)$. Then $g_{0*} = g_{1*}$
        and $g_{0}^* = g_1^*$.\end{theorem}

\begin{lemma}
    \label{lem:Homotopy equiv Homotopy 2}
	Theorem \ref{thm:Homotopy 2} is equivalent to Theorem \ref{thm:Homotopy}.
\end{lemma}
\begin{proof}
	Theorem \ref{thm:Homotopy} $\implies$ Theorem \ref{thm:Homotopy
        2}$\colon g_0$ and $g_1$ are homotopic by the homotopy
        $H:(X,\mathcal{F}_X) \times (I,\tau) \to (X,\mathcal{F}_X) \times
        (I,\tau)$ defined by $H(x,t) = (x,t)$, which is uniformly continuous, since the
        identity from a semi-uniform space to itself is uniformly continuous. 
.
	
	Theorem \ref{thm:Homotopy 2} $\implies$ Theorem
        \ref{thm:Homotopy}$\colon$ Suppose $H:(X,\mathcal{F}_X) \times (I,\tau)
        \to (Y,\mathcal{F}_Y)$ is a homotopy 
        for $f_0, f_1:(X,\mathcal{F}_X) \to (Y,\mathcal{F}_Y)$. Then $f_i = H
        \circ g_i$. By functoriality, we have
        $f_{i*} = H_*g_{i*}$ and $f^*_{1} = g^*_i H^*$. Since $g_{0*} = g_{1*}$
        and $g^*_0 = g^*_1$, by hypothesis, the result follows. 
\end{proof}

\subsubsection{Acyclic carriers}
Before proceeding with the proof of Theorem \ref{thm:Homotopy 2}, we will first
review some terminology and results from\cite{Eilenberg_Steenrod_1952}, VI, section 5 on common acyclic carriers in simpicial
complexes that we need for the proof. More details may be found there.

\begin{definition}
    Let $K$ be a simplicial complex and $G$ an abelian group, let
    $C_0(K;G)$ denote the $0$-th dimensional simplicial chain complex,
    generated by the vertices of $K$, and let $c \coloneqq \sum_i g_i v_i \in
    C_0(K;G)$ be an arbitrary element of $C_0(K;G)$. We define the index map $\Ind:C_0(K;G) \to
    G$ by $\Ind(c) = \Ind\left(\sum_i g_i v_i\right) \coloneqq \sum_i g_i$. 
\end{definition}

\begin{definition}
We say that a map $f:K \to K'$ between simplicial complexes is a \emph{chain
map} if the induced map $f_*:C_*(K) \to C_*(K')$ satisfies $df = fd$. We say
that $f$ is an \emph{algebraic map} if $f$ is a chain map and $Ind(f(c)) = Ind(c)$ for all $c \in C_0(K)$.
\end{definition}

\begin{definition}
A function which maps simplices $\sigma$ of a simplicial complex $\Sigma$ to
subcomplexes $C(\sigma) \subset \Sigma'$ of a simplicial complex $\Sigma'$ is
called a \emph{carrier function} if, for every face $\gamma \subset \sigma$,
$C(\gamma)$ is a subcomplex of $C(\sigma)$. If $f:\Sigma \to \Sigma'$ is an
algebraic map such that $\alpha \subset
\sigma \in \Sigma$ implies $f(\alpha) \subset C(\sigma)$, then $C$ is called the \emph{carrier of $f$}. 
\end{definition}

\begin{definition}
    We say that a simplicial complex $K$ is \emph{acyclic} if $\tilde{H}_q(K) =
    0$, where $\tilde{H}_*(K)$ is the reduced simplicial homology of the
    simplicial complex $K$. We say that a carrier function $C$ is
    \emph{acyclic} if, for each simplex $\sigma$, the simplicial complex
    $C(\sigma)$ is acyclic. 
\end{definition}

\begin{proposition}[\cite{Eilenberg_Steenrod_1952}, Theorem 5.8]
\label{prop:Common acyclic carrier}
Suppose $K$ and $K'$ be simplicial complexes and suppose that $f,g:K \to K'$ are
algebraic maps with a common acyclic carrier function $C$. Then $f_*=g_*$ and $f^* = g^*$.
\end{proposition}

\subsubsection{Proof of the homotopy axiom}

In this section, we will give the proof of Theorem \ref{thm:Homotopy 2}, which
will follow from the next three lemmas.

\begin{lemma}
    \label{lem:VR complex of interval}
    Given $r>0$, let $U_r \subset I \times I$ be defined as in Example
    \ref{ex:Semi-uniform examples}.(\ref{ex:Semi-uniformity from metric}). For any $r>0$, the Vietoris-Rips complex of $(I,U_r)$, $\Sigma^I_{U_r}$, is acyclic.
\end{lemma}

\begin{proof}
This follows directly from \cite{Hausmann_1995}, Lemma 2.2.
\end{proof} 

\begin{lemma}
	\label{lem:The critical acyclic carrier}
Let $X$ be a set with relation $R \subset X \times X$. Consider $I$ with the
Euclidean metric and the relation $U_r \subset I \times I$ defined above, and
let $\sigma = \{v_0,\dots,v_k\}$ be a simplex of $\Sigma^X_R$. Suppose
$S(\sigma)$ is the subcomplex of $\Sigma^{X\times I}_{R \times U_r}$ which consists
of all of the simplices in $\Sigma^{X \times I}_{R \times U_r}$ whose vertices are of the 
form $(v_i,t)$, $v_i \in \sigma$, $t\in I$. Then $S(\sigma)$ is acyclic.
\end{lemma}
\begin{proof}Consider a simplex $\gamma= ((v_0,t_0),\dots,(v_n,t_n))$ in
    $S(\sigma)$.
    By hypothesis, $\sigma \in \Sigma^X_R$, so $\gamma$ is a simplex in $S$ iff $(t_0,\dots,t_n)\in
    \Sigma^I_{U_r}$. Since $\sigma \in \Sigma^X_R$, however, this implies
    that the map $\tilde{f}:S(\sigma) \to S(\sigma)$ induced by the map $(v_i,t) \mapsto
    (v_0,t)$ is contiguous with the identity, and therefore the homology and
    cohomology of $S(\sigma)$ are equal to that of $I \times \{v_0\} \cong I$. Since $I$ is
    acyclic by Lemma \ref{lem:VR complex of interval}, the result follows. 
\end{proof}

\begin{lemma}
    \label{lem:Homotopy invariance for one controlled set}
    Let $(X,\mathcal{F}_X)$ be a semi-uniform space, and consider the continuous maps
    $g_0,g_1: (X,\mathcal{F}_X) \to (X,\mathcal{F}_X) \times (I,\tau)$ given by
\begin{align*}
    g_0(x) &\coloneqq (x,0),\\
    g_1(x) &\coloneqq (x,1),
\end{align*}
where $\tau = \mathcal{F}_0$ is the topological semi-uniform structure. Then, for every $V \in \tau$ and $U \in \mathcal{F}_X$, there exists a $\tilde{U}
\in \mathcal{F}_X$ be such that $g_i(\tilde{U}) \subset U \times V$ for
$i \in \{0,1\}$, and 
\begin{align*}
    g_{0\tilde{U}*} &= g_{1\tilde{U}*}\colon H_*^{VR}(X\colon \tilde{U}) \to H_*^{VR}(X \times I\colon U \times V)\\
g_{0\tilde{U}}^*&= g_{1\tilde{U}}^*\colon H_{VR}^*(X\times I\colon U \times V) \to H_{VR}^*(X
    \colon \tilde{U})
\end{align*}
\end{lemma}

\begin{proof}
    Since each $g_i$, $i \in \{0,1\}$ is continuous, there exist $U_i \in
    \mathcal{F}_X$ such that
    $g_i(U_i) \subset V$. Define $\tilde{U} = U_0 \cap U_1$, and note that
    $\tilde{U} \in \mathcal{F}_X$ by definition of a semi-uniform structure. Now, for each
    simplex $\sigma$ of $\Sigma^X_{\tilde{U}}$, consider the subcomplex
    $C(\sigma)$ of $\Sigma^{X\times I}_{U \times V}$ consisting of all
    simplices with vertices of the form $(x_i,t)$, where each $x_i$ is a vertex of
    $\sigma$ and $t \in I$. By Lemma \ref{lem:The critical acyclic carrier},
    $C(\sigma)$ is acyclic. Note that, if $\gamma$ is a face of $\sigma$, then $C(\gamma)
    \subset C(\sigma)$, the $g_i$ are algebraic maps, since $g_{i_*}(\alpha
    x_i) = \alpha g_i(x_i)$, so the index of any $0$-chain is preserved, and,
    furthermore, $\alpha \subset \sigma \in \Sigma^X_{\tilde{U}}$ implies that
    $f(\alpha) \subset C(\sigma)$. Since $\sigma \in
    \Sigma^X_{\tilde{U}}$ is arbitrary, we have that $C$ is
    a common acyclic carrier for the maps $\tilde{g}_0$ and $\tilde{g}_1$, and
    the conclusion follows from Proposition \ref{prop:Common acyclic carrier}. 
\end{proof}

We now return to the proofs of Theorem \ref{thm:Homotopy 2}, which implies
Theorem \ref{thm:Homotopy}.

\begin{proof}[Proof of Theorem \ref{thm:Homotopy 2}]
By construction of the product semi-uniform structure, sets of the form $U
\times V$, $U \in \mathcal{F}_X, V\in \tau$ generate $\mathcal{F}_X \times
\tau$, and therefore Theorem \ref{thm:Homotopy 2} follows from Lemma
\ref{lem:Homotopy invariance for one controlled set} by taking the appropriate
limits. 
\end{proof}

\begin{proof}[Proof of Theorem \ref{thm:Homotopy}] Theorem \ref{thm:Homotopy} now follows from Lemma \ref{lem:Homotopy
    equiv Homotopy 2} and Theorem \ref{thm:Homotopy 2}.
\end{proof}

\section{Isomorphism theorems}
\label{sec:Isomorphisms}

In this section, we give a number of results which aid in the calculation of
the Vietoris-Rips homology and cohomology groups. We begin by giving a simple alternative proof of Theorem
5.1 in \cite{Hausmann_1995}, using results in
\cite{Dowker_1952}. 
\begin{theorem}
    For the category of compact metric pairs (i.e. pairs of metric spaces $(X,A)$ with $X$ and
    $A$ compact and $A \subset X$), $H^*_{VR}(X,A,\mathcal{F}_0) \cong
    \check{H}^*(X,A)$, where $\check{H}^*(X,A)$ is the \v Cech cohomology of the pair
    $(X,A)$.
\end{theorem}

\begin{proof}
    By Proposition \ref{prop:Metric and topological structures on compact
    spaces}, $\mathcal{F}_0 = \mathcal{F}_\tau$, and therefore
    $H^*_{VR}(X,A,\mathcal{F}_0)
    \cong H^*_{VR}(X,A,\mathcal{F}_\tau)$. However, the groups
    $H^*_{VR}(X,A,\mathcal{F}_\tau)$ are the Alexander cohomology groups from
    \cite{Dowker_1952}. By Theorem 2 in \cite{Dowker_1952}, then, we have
    $H^*_{VR}(X,A,\mathcal{F}_\tau) \cong \check{H}^*(X,A)$.
\end{proof}

We continue with several situations where the semi-uniform structure has a
maximal element, facilitating the calculation of the Vietoris-Rips homology and
cohomology groups.

\begin{thm}
	\label{thm:H_VR of finite metric space}
	Let $(X,d)$ be a finite semi-pseudometric space and suppose that $q>0$ is a
        positive real number.  Let
        $U_{\leq q} = \{ (x,y) \mid d(x,y) \leq q \}$. Then
        \begin{align*} H^{VR}_*(X,\mathcal{F}_q) &\cong H_*\mleft(\Sigma^X_{U_{\leq
                q}}\mright), \text{
            and}\\
            H_{VR}^*(X,\mathcal{F}_q) &\cong H^*\mleft(\Sigma^X_{U_{\leq
                q}}\mright).
        \end{align*}
\end{thm}
\begin{proof}
    For every $\epsilon > 0$, define 
    \[ V_{q+\epsilon} \coloneqq \{ (x,y) \in X \times X \mid d(x,y) < q +
    \epsilon\}. \]
    Since $X$ is finite, there is an $\epsilon_0 > 0$ such that, for any $0 <
    \epsilon < \epsilon_0$, $V_{q+\epsilon} = U_{\leq q}$. Since the collection
    of sets $\{V_{q + \epsilon} \mid 0< \epsilon <
    \epsilon_0\}$ is cofinal in $\mathcal{F}_q$ and, for each such
    $0<\epsilon<\epsilon_0$, we have 
    \begin{align*} 
        H^{VR}_*(X,V_{q+\epsilon}) &\cong H_*\mleft(\Sigma^X_{U_{\leq q}}\mright), \text{ and}\\
        H_{VR}^*(X,V_{q+\epsilon}) &\cong H^*\mleft(\Sigma^X_{U_{\leqq}}\mright).
    \end{align*}
    the result follows.
\end{proof}

\begin{remark}
    Recall that, for the classical $(<)$- and $(\leq)$-Vietoris-Rips complexes
    $\Sigma_{VR<}(X)$ and $\Sigma_{VR\leq}(X)$ defined at a
    scale $q > 0$, the set $\{x_0,\dots,x_q\}$ is a simplex in $\Sigma_{VR}(X)$
    iff, for every pair $(x_i,x_j)$, $d(x_i,x_j) < q$ or $d(x_i,x_j) \leq q$,
    respectively. The above theorem shows that, for finite pseudo-metric
    spaces, the
    classical Vietoris-Rips homology and the one defined here are the
    isomorphic for
    the $\leq$ version of the Vietoris-Rips complex, and also for the
    $<$-Vietoris-Rips complexi for all $q>0$ such that $d(x,x') \neq q$
    for any $x,x' \in X$.
\end{remark}

\begin{thm}
    Let $M$ be a closed Riemannian manifold, and let $d_{GW}$ denote the Gromov-Hausdorff
    distance. There exists an $\epsilon_0>0$ such that, for any
    $0<\epsilon<\epsilon_0$, there exists a $\delta>0$ where, if $(X,d)$ is a
    finite metric space with $d_{GW}(M,X)< \delta$, then
    $H_*^{VR}(X,\mathcal{F}_\epsilon) \cong H_*(M)$ and $H^*_{VR}(X,\mathcal{F}_\epsilon) \cong H^*(M)$.
\end{thm}

\begin{proof}
    
    By Theorem 1.1 in \cite{Latschev_2001} there exists an $\epsilon_0>0$
    such that for any $0<\epsilon
    < \epsilon_0$, $H_*^{VR}(X,U_\epsilon) \cong H_*(M)$ and
    $H^*_{VR}(X,U_\epsilon) \cong H^*(M)$, where $U_\epsilon$ is defined by
    \[ U_\epsilon \coloneqq \{ (x,y) \in X \times X \mid d(x,y) <
    \epsilon\}.\] 

    Now fix an $0 < \epsilon < \epsilon_0$, and let $0< \epsilon'< \epsilon_0 -
    \epsilon$. Then $H_*^{VR}(X,U_{\epsilon+\epsilon'}) \cong H_*(M)$ and
    $H^*_{VR}(X,U_{\epsilon+\epsilon'}) \cong H^*(M)$. Since the
    $U_{\epsilon+\epsilon'}$ are cofinal in $\mathcal{F}_\epsilon$, the result
    follows.
\end{proof}

\begin{thm}
    Let $G = (V,E)$ be an undirected graph, and let $(V,\mathcal{F}_E)$ be the
    associated semi-uniform space, i.e. $\mathcal{F}_E$ is the filter generated
    by $E \in V \times V$. Then $H_*^{VR}(V,c_E)$ and $H^*_{VR}(V,c_E)$ are equal to the homology and cohomology, respectively, of the clique complex of $G$.
\end{thm}
\begin{proof}
    By construction, $E$ is a maximal element of $\mathcal{F}_E$ and therefore
    the single-element family $\{E\}$ is cofinal in $\mathcal{F}_E$, which
    implies the theorem.
\end{proof}

\acknowledgement
We are grateful to Luis Jorge Palacios Vela and Henry Adams for interesting discussions and comments.

\bibliography{/home/antonio/Bib/all.bib}

\end{document}